\documentclass[11pt]{amsart}
\usepackage[pagewise]{lineno}
\usepackage{graphicx}
\usepackage{amsmath,amsthm,amssymb}
\usepackage{amsfonts}
\usepackage{indentfirst}

\newtheorem{thm}{Theorem}[section]

\theoremstyle{definition}
\newtheorem{defn}{Definition}[section]

\theoremstyle{remark}

\begin{document}

	\title[Pathological Foliations]{A Simple Example of  Pathological Foliations in Skew-Product Diffeomorphisms}
	\author{Zhihong Xia, Peizheng Yu}
	\address{Department of Mathematics, Northwestern University, Evanston, 
		IL 60208 USA}
	\address{School of Mathematics, Shandong University, Jinan, China}
	\email{xia@math.northwestern.edu, yupzh@sdu.edu.cn}
	\date{version2: October 31, 2024}
	\maketitle

	\begin{abstract}
		Inspired by examples of Katok and Milnor \cite{Milnor1997}, we
		construct simple examples of skew-product volume preserving
		diffeomorphism where the center foliation is pathological in the
		sense that, there is a full measure set whose intersection with any
		center leaf contains at most one point. Comparing with other
		examples in literature, our mechanism for producing such
		pathological foliation is simple and elementary.
	\end{abstract}

\section{Introduction}

An interesting and intriguing phenomenon in dynamical systems is the
pathological foliations. Roughly speaking, a foliation is pathological
if there is a full volume set that meets every leaf of the foliation
on a set of leaf-volume zero. In fact, there are examples with a full
measure set that intersects each leaf only in a finite number of
points. In \cite{Milnor1997}, J. Milnor constructed an example,
inspired by A. Katok, of such a non-absolutely continuous foliation on
the unit square. In \cite{ShubWilkinson2000}, M. Shub and A. Wilkinson
found the same phenomenon on $\mathbb{T}^3$ with a different
approach. They termed this phenomenon aptly {\em ``Fubini's
	nightmare''}\/. Surprisingly, this phenomenon is persistant and
robust under perturbations. We refer readers to a survey
\cite{Ures2007} written by F. R. Hertz, M. R. Hertz and R. Ures for
more details. Saghin and Xia \cite{RaduXia2009}, using a different
mechanism, showed more examples of systems with persistent
non-absolute continuous center and weak unstable foliations, where
these foliations are not necessarily compact. Pesin \cite{Pesin2004}
also showed some examples of a non-absolutely continuous foliation.

In this paper, we construct a simple example of skew-product diffeomorphism where the center foliation is pathological. And our main theorem is as follows:
\begin{thm}
	There exists a full measure set $E$ on $[0,1] \times \mathbb{T}^2$, together with a family of disjoint curves $\Gamma_\beta$ which fill out $[0,1] \times \mathbb{T}^2$, so that each curve $\Gamma_\beta$ intersects the set $E$ at most a single point.
\end{thm}

The idea of the proof of this theorem is inspired by Milnor's proof of
Katok's paradoxical example, see \cite{Milnor1997}. We construct a
path of Anosov area-preserving diffeomorphism $f_p$, $p \in [0,1]$ on
$\mathbb{T}^2$ beginning with Arnold's cat map. One can certainly
replace the cat map with any Anosov linear automorphisms on
$\mathbb{T}^2$. We construct Markov partitions for each $f_p$ such
that each Markov partition has two rectangles and the Lebesgue measure
of one rectangle is strictly increasing with respect to $p$. Then we
can construct a full measure set $E$ and a family of disjoint curves
$\Gamma_\beta$ satisfying the conditions in the theorem by using a
(semi)conjugacy from a shift space.

F. R. Hertz, M. R. Hertz and R. Ures \cite{Ures2007} showed similar
results by a different mechanism. They perturbed eigenvalues of fixed
points, instead of volume of Markov partitions. Their results are
heavily based on Lyapunov exponents and other techniques in ergodic
theory, as are all results in pathological foliations in literature
(cf.\ Shub-Wilkinson \cite{ShubWilkinson2000}, Saghin-Xia
\cite{RaduXia2009}, Pesin \cite{Pesin2004}). Our construction is
direct and elementary, no Lyapunov exponents are involved.

In section 2, we discuss about some definitions and properties of
Anosov diffeomorphisms and Markov partitions and especially introduce
Arnold's cat map. In section 3, we prove our main theorem. The main
idea is to perturb the Arnold's cat map in such a way that the areas
of Markov petitions of the perturbed map are different.

\section{Preliminaries}
Our proof of the main theorem uses some basic properties of a hyperbolic toral automorphism which is an Anosov diffeomorphism. We will introduce some basic definitions in this section. The content of this section is based on Brin and Stuck's book \cite{BrinStuck2002}.

\begin{defn}
	Let $M$ be a smooth Riemannian manifold and $f:M\rightarrow M$ be a diffeomorphism. A compact, $f$-invariant subset $\Lambda\subset M$ is called {\em hyperbolic}\/ if there are $\lambda\in(0,1),C>0,$ and families of subspaces $E^s(x)\subset{T_xM}$ and $E^u(x)\subset{T_xM}$, $x\in \Lambda$, such that for every $x\in \Lambda$,
	\begin{enumerate}
		\item $T_xM=E^s(x)\oplus E^u(x),$
		\item $\|df_x^nv^s\|\leq C\lambda^n\|v^s\|\text{ for every }v^s\in E^s(x)\mathrm{~and~}n\geq0,$
		\item $\|df_x^{-n}v^u\|\leq C\lambda^n\|v^u\|\text{ for every }v^u\in E^u(x)\mathrm{~and~}n\geq0,$
		\item $df_{x}E^{s}(x)=E^{s}(f(x))\text{ and }df_{x}E^{u}(x)=E^{u}(f(x)).$
	\end{enumerate}
	In particular, if $\Lambda=M$, then $f$ is called an {\em Anosov diffeomorphism}\/
\end{defn}

Hyperbolic toral automorphisms are examples of Anosov diffeomorphisms. Explicitly,
\begin{defn}
	Let $M=\mathbb{T}^{n}$. Consider an $n \times n$ matrix $A$ with determinant one and with integer entries. The matrix $A$ induces a {\em toral automorphism}\/: $f_A:\mathbb{T}^n=\mathbb{R}^n/\mathbb{Z}^n\to\mathbb{T}^n$ defined by $f_{A}(x)=Ax\bmod\mathbb{Z}^{n}$. 
	
	Moreover, if all eigenvalues of $A$ are away from the unit circle, then $f_{A}$ is a {\em hyperbolic toral automorphism}\/. 
\end{defn}

\

The best known example of a hyperbolic toral automorphism is Arnold's cat map that is $f_A:\mathbb{T}^2\to\mathbb{T}^2$ where
\begin{equation*}
	A=\begin{pmatrix}2&1\\1&1\end{pmatrix}. 
\end{equation*}
And the eigenvalues of $A$ is $\lambda=\frac{3+\sqrt{5}}{2}>1$ and $\lambda^{-1}=\frac{3-\sqrt{5}}{2}$. The corresponding eigenvectors are $\nu_{\lambda}=(1,\frac{\sqrt{5}-1}{2})$ and $\nu_{\lambda^{-1}}=(1,\frac{-\sqrt{5}-1}{2})$.

To figure out the features of the Arnold's cat map, we need to introduce more concepts.
\begin{defn}
	Let $f:M\rightarrow M$ be an Anosov diffeomorphism. For every $x\in M$, the {\em (global) stable}\/ and
	{\em unstable manifolds}\/ of $x$ are defined by
	\begin{equation*}
		\begin{aligned}W^s(x)&:=\{y\in M\colon \mathrm{dist}(f^n(x),f^n(y))\to0\text{ as }n\to\infty\},\\W^u(x)&:=\{y\in M\colon \mathrm{dist}(f^{-n}(x),f^{-n}(y))\to0\text{ as }n\to\infty\}.\end{aligned}
	\end{equation*}
	And for any $\rho>0$, we define
	\begin{equation*}
		\begin{aligned}W^s(x,\rho)&:=\{y\in M:\mathrm{dist}(f^n(x),f^n(y))<\rho,\forall n\in\mathbb{N}_0\},\\W^u(x,\rho)&:=\{y\in M:\mathrm{dist}(f^{-n}(x),f^{-n}(y))<\rho,\forall n\in\mathbb{N}_0\}.\end{aligned}
	\end{equation*}
\end{defn}

With these definitions, we can check that if $x_0$ is a fixed point of $f$, we have
\begin{equation*}
	\begin{aligned}
		&W^s(x_0)=\bigcup_{n\in\mathbb{N}}\left(f\right)^{-n}\left(W^s(x_0,\rho)\right)\\
		&W^u(x_0)=\bigcup_{n\in\mathbb{N}}\left(f\right)^n\left(W^u(x_0,\rho)\right).
	\end{aligned}
\end{equation*}

\
An important property of an Anosov diffeomorphism is called structurally stability.
\begin{defn}
	Let $M$ be a smooth Riemannian manifold. A $C^1$ map $f:M\rightarrow M$ is called {\em structurally stable}\/ if there exists a neighborhood $U$ of $f$ in the $C^1$ topology  such that every $g \in U$ is topologically conjugate to $f$.
\end{defn}

\begin{thm}
	Anosov diffeomorphisms are structurally stable.
\end{thm}

Another concept that we would like to introduce is called a Markov partition. It allows us to study the dynamics of $f$ using symbolic dynamics.
\begin{defn}
	Let $f:M\rightarrow M$ be an Anosov diffeomorphism. A collection of subset of $M$, $\mathcal{R}=\{{R_{1}},\ldots,{R_{n}}\}$ is called a {\em Markov partition} for $(M,f)$ if
	\begin{enumerate}
		\item $\operatorname{cl}(\operatorname{int}R_i)=R_i\text{ for each }R_i;$
		\item $\operatorname{int}R_i\cap\operatorname{int}R_j=\varnothing\operatorname{for}i\neq j;$
		\item $M=\bigcup_i R_i;$
		\item If $f^m(\operatorname{int}R_i)\cap \operatorname{int}R_j \neq \varnothing$ for some $m\in \mathbb{Z}$ and $f^n(\operatorname{int}R_j)\cap \operatorname{int}R_k \neq \varnothing$ for some $n\in \mathbb{Z}$, then $f^{m+n}(\operatorname{int}R_i)\cap \operatorname{int}R_k \neq \varnothing$.	
	\end{enumerate}
\end{defn}

As an example, we can construct a Markov partition for Arnold's cat map $f_A$ by draw segments of stable and unstable manifolds of the fixed point $(0,0)$ until they cross sufficiently many times and separate $\mathbb{T}^2$ into two disjoint rectangles $R_0$ and $R_1$: $R_0$ consists of two parts $B_1$ and $B_2$; $R_1$ consists of three parts $A_1$, $A_2$ and $A_3$. See figure 1.

\begin{figure}[htbp]
	\centering
	\includegraphics[width=0.7\textwidth,height=0.415\textheight]{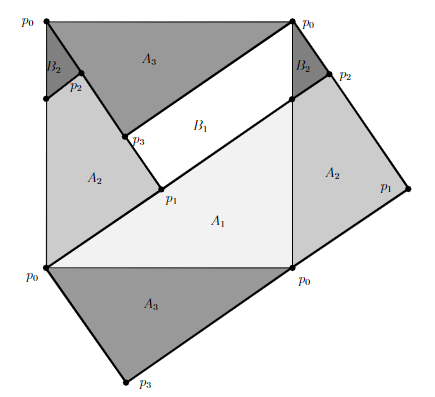}
	\caption{Markov partition of Arnold's cat map}
	\label{}
\end{figure}

If we consider the fundamental domain in Figure 1 as $[0,1]\times[0,1]$, the coordinates of $p_1$, $p_2$, $p_3$ in $[0,1]\times[0,1]$ are as follows:
\begin{equation*}
	\begin{aligned}
		&p_1=(\frac{1}{\sqrt{5}},\frac{\sqrt{5}-1}{2\sqrt{5}})\\
		&p_2=(\frac{3-\sqrt{5}}{2\sqrt{5}},1+\frac{-\sqrt{5}-1}{2}\cdot\frac{3-\sqrt{5}}{2\sqrt{5}})\\
		&p_3=(\frac{2}{5+\sqrt{5}},1+\frac{-\sqrt{5}-1}{2}\cdot\frac{2}{5+\sqrt{5}})
	\end{aligned}
\end{equation*}
Since $(0,0)$,$(1,0)$,$(0,1)$,$(1,1)$ in $[0,1]\times[0,1]$ represent the same fixed point on $\mathbb{T}^2$, we denote them by $p_0$. And in order to distinguish them in $[0,1]\times[0,1]$, we write $p_0(i,j)$ to be the point $(i,j)$, $i,j=0,1$. 

Note that the preimage of the segment $\overline{p_0(1,1)p_3}=\{(x,y)\big| y-1=\frac{\sqrt{5}-1}{2}(x-1), x\in [\frac{2}{5+\sqrt{5}},1]\}$ is contained in the segment $\{(x,y)\in \overline{p_0(1,1)p_3} \big| x\in [1-\frac{2}{5+\sqrt{5}},1]\}$. 

The preimage of the segment $\overline{p_0(0,0)p_2}=\{(x,y)\big| y=\frac{\sqrt{5}-1}{2}x, x\in [0,1+\frac{3-\sqrt{5}}{2\sqrt{5}}]\}$ is contained in the segment $\{(x,y)\in\overline{p_0(0,0)p_2} \big| x\in [0,1-\frac{2}{5+\sqrt{5}}]\}$. 

And the image of the segment $\overline{p_0(0,1)p_1}=\{(x,y)\big| y-1=\frac{-\sqrt{5}-1}{2}x, x\in [0,\frac{1}{\sqrt{5}}]\}$ is contained in $\overline{p_0(0,1)p_1}$ where the $x$-coordinate of $p_1$ is less than $1-\frac{2}{5+\sqrt{5}}$. 

Thus, if we perturb Arnold's cat map into $f_p$ in the region $\{(x,y)\big|x\in (a,b)\}$, where $1-\frac{2}{5+\sqrt{5}}<a<b<1$. Then $(0,0)$ is still a fixed point of $f_p$ and if we draw stable manifold and unstable manifold of $(0,0)$, we can get $\overline{p_0(0,0)p_2}$ and $\overline{p_0(0,1)p_1}$ will not change. And $\overline{p_0(1,1)p_3}$ will be perturbed for the perturbed map $f_p$. This inspires us the proof of our main theorem.

\section{Proof of Main theorem}
\textbf{Step 1:} Our first step is to construct a continuous path of area-preserving Anosov diffeomorphism $f_p$, $p \in [0,1]$ on $\mathbb{T}^2$ beginning with the Arnold's cat map.

\

Let $f_A$ be the Arnold's cat map on $\mathbb{T}^2$ that is
\begin{equation*}
	f_A(x,y)=(2x+y,x+y)\mathrm{~mod~}1.
\end{equation*}

Take real number $1-\frac{2}{5+\sqrt{5}}<a<1$ and $\delta>0$ small enough such that $1-\frac{2}{5+\sqrt{5}} <a-\delta <a+\delta<1$. Then we fix a smooth bump function $C(x)$ on $[0,1]$ such that
\begin{enumerate}
	\item $C(x)\equiv0$ on $[0,a-\delta]\cup [a+\delta,1]$;
	
	\item $0<C(x)\leq1$ on $(a-\delta,a+\delta)$.
\end{enumerate}

Now, for any $p \in [0,1]$, we define a function $\phi_p(x)$ on $[0,1]$ such that $\phi_p(x)=p\epsilon_0C(x)$. Then we define a diffeomorphism $f_p$ on $\mathbb{T}^2=\mathbb{R}^2/\mathbb{Z}^2$ such that
\begin{equation*}
	f_p(x,y)=(2x+y-\phi_p(x),x+y-\phi_p(x))\mathrm{~mod~}1.
\end{equation*}
Note that $f_0=f_A$. And for a fixed sufficiently small $\epsilon_0>0$, $f_p$ is $C^1$-closed to $f_A$ for any $p \in [0,1]$. Thus by structurally stability of Anosov diffeomorphisms, we can get $f_p$ is an Anosov diffeomorphism for any $p \in [0,1]$. Also, for any $p \in [0,1]$, it follows from the Jacobian of $f_p$ equal to 1 that $f_p$ is an area-preserving diffeomorphism. Moreover, $f_p(x,y)=f_A(x.y)$ for $x(\mathrm{~mod~}1) \in [0,a-\delta]\cup [a+\delta,1]$. In conclusion, $f_p$, $p \in [0,1]$ is a path of area-preserving Anosov diffeomorphism.

\

\textbf{Step 2:} Our second step is to construct Markov partitions $\mathcal{R}_p$ with two rectangles for $f_p$ for each $p \in [0,1]$.

\

For any $p \in [0,1]$, since $(0,0)$ is a fixed point of $f_p$, one way to construct a Markov partition for $f_p$ is to draw segments of $W_p^u((0,0))$ and $W_p^s((0,0))$ until they cross sufficiently many times and separate $\mathbb{T}^2$ into two disjoint (curvilinear) rectangles. Explicitly, we view $\mathbb{T}^2 = [0,1]\times[0,1]/\sim$, where we identify $(0,y)\sim (1,y)$ and $(x,0)\sim (x,1)$. Then we construct Markov partitions as follows (see Figure 2):
\begin{enumerate}
	\item From $(0,0)$ draw the line $l_p(0,0)=W_p^u((0,0))$ in the unit square and stops when it hits the boundary of the unit square.
	
	\item Form $(0,1)$ draw the line $l_p(0,1)=W_p^s((0,0))$ in the unit square and stops when it hits $l_p(0,0)$.
	
	\item From $(1,1)$ draw the line $l_p(1,1)=W_p^u((0,0))$ in the unit square and stops when it hits $l_p(0,1)$.
	
	\item Finally using symmetry, draw the extension $l_p$ of the line $l_p(0,0)$ in the unit square and stops when it hits $l_p(0,1)$ to complete the rectangles.
\end{enumerate}

\begin{figure}[htbp]
	\centering
	\includegraphics[width=0.65\textwidth,height=0.38\textheight]{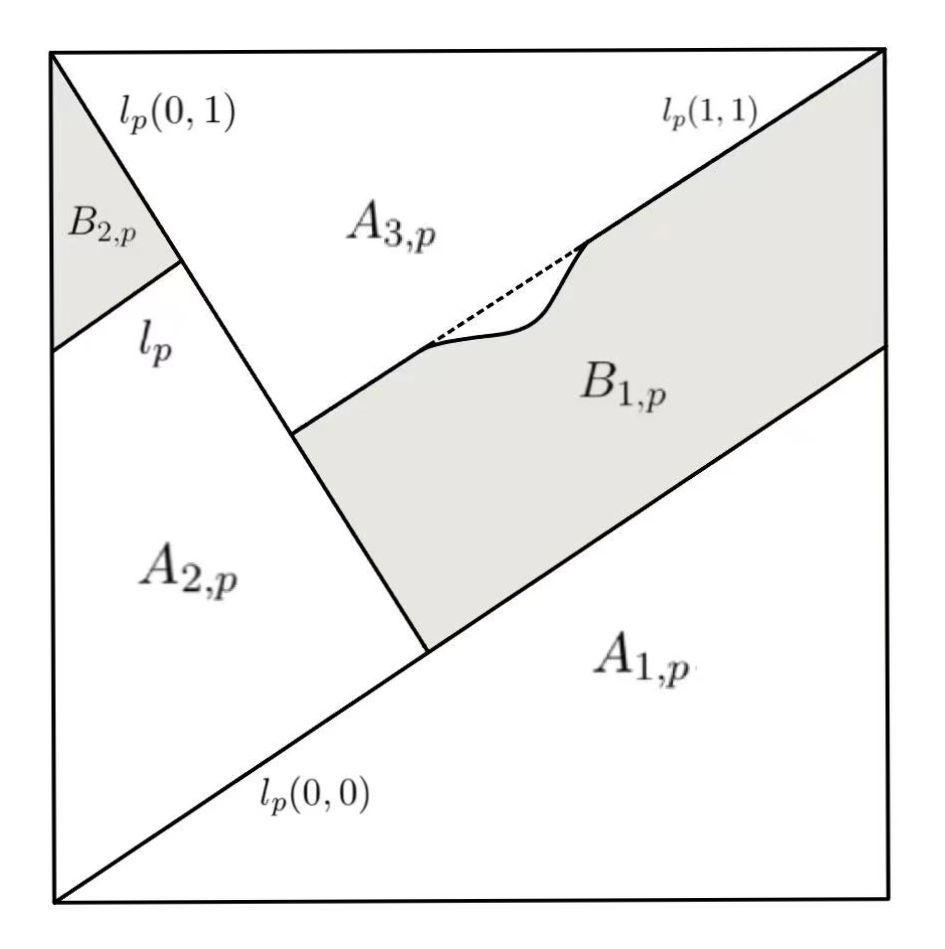}
	\caption{Markov partition of $f_p$}
	\label{}
\end{figure}

We denote the rectangle consisting of $B_{1,p}$ and $B_{2,p}$ by $R_p^0$, and denote the rectangle consisting of $A_{1,p}$, $A_{2,p}$ and $A_{3,p}$ by $R_p^1$. Then $\mathcal{R}_p=\{R_p^0, R_p^1\}$ is a Markov partition for $f_p$.

Note that when $p=0$, the Markov partition $\mathcal{R}_0$ is exactly the Markov partition for the Arnold's cat map we introduced in last section. And since for any $p\in [0,1]$,
\begin{equation*}
	\begin{aligned}
		&W_p^s((0,0))=\bigcup_{n\in\mathbb{N}}\left(f_p\right)^{-n}\left(W_p^s((0,0),\rho)\right)\\
		&W_p^u((0,0))=\bigcup_{n\in\mathbb{N}}\left(f_p\right)^n\left(W_p^u((0,0),\rho)\right).
	\end{aligned}
\end{equation*}
and for a sufficient small $\rho >0$, $W_p^i((0,0),\rho)=W_0^i((0,0),\rho)$, $i=u,s$, we can find how the stable manifold and unstable manifold of $(0,0)$ change when $p$ changes from 0 to 1.

Recall that $\forall p \in [0,1]$,  $f_p(x,y)=f_0(x,y)$ for $x(\mathrm{~mod~}1) \in [0,a-\delta]\cup [a+\delta,1]$. By our construction of $f_p$, since $a-\delta > 1- \frac{2}{5+\sqrt{5}}$, we can check that $l_p(0,0)=l_0(0,0)$, $l_p(0,1)=l_0(0,1)$ and $l_p=l_0$. And next let's figure out how $l_p(1,1)$ changes when $p$ changes.

To simplify the calculation, let's consider the fundamental domain as $[-1,0]\times[-1,0]$. Then $l_p(1,1)$ in this fundamental domain is a line from $(0,0)$ and stops when it hits the stable manifold of $(0,0)$. Let $\lambda=\frac{3+\sqrt{5}}{2}$, then we can see that $l_0(1,1)$ in this fundamental domain is the segment $\{(x,y)\big| y=\frac{-1+\sqrt{5}}{2} x, x\in[-\lambda \frac{2}{5+\sqrt{5}},0] \}$. Thus, by our construction of $f_p$, for any $p \in[0,1]$, the segment 
\begin{equation*}
	\begin{aligned}
		l_p(1,1)= &f_p(\{(x,y)\big| y=\frac{-1+\sqrt{5}}{2} x, x\in[-\frac{2}{5+\sqrt{5}},0] \})\\
		=&\begin{cases}
			l_0(1,1)\quad\text{  for }x\in [-\lambda \frac{2}{5+\sqrt{5}},\lambda(a-\delta-1)]\cup[\lambda(a+\delta-1),0]\\
			f_p(\{(x,y)\big| y=\frac{-1+\sqrt{5}}{2} x, x\in[a-\delta-1,a+\delta-1] \}),\text{ otherwise}
		\end{cases}
	\end{aligned}
\end{equation*}

Moreover, if we write $l_p(1,1)=\{(X,Y(X))\big|X\in[-\lambda \frac{2}{5+\sqrt{5}},0]\}$, we can calculate the area of the domain bounded by $l_p(1,1)$ and $l_0(1,1)$ for each $p \in [0,1]$: 
\begin{equation*}
	\begin{aligned}
		\text{Area}_p&=\int_{a-\delta-1}^{a+\delta-1}\frac{\sqrt5+1}{2}\frac{\sqrt5+3}{2}xdx-\int_{\lambda(a-\delta-1)}^{\lambda(a+\delta-1)}Y_p(X)dX\\
		&=\int_{a-\delta-1}^{a+\delta-1}\frac{\sqrt5+1}{2}\frac{\sqrt5+3}{2}xdx-\int_{a-\delta-1}^{a+\delta-1}(\frac{\sqrt5+1}{2}x -\phi_p(x))(\frac{\sqrt5+3}{2}-\frac{d\phi_p}{dx}(x))dx\\
		&=\int_{a-\delta-1}^{a+\delta-1}\phi_p(x)dx+\int_{a-\delta-1}^{a+\delta-1}\frac{\sqrt5+1}{2}(x\phi_p(x))'dx-\int_{a-\delta-1}^{a+\delta-1}\frac{1}{2}(\phi_p^2(x))'dx\\
		&=\int_{a-\delta-1}^{a+\delta-1}\phi_p(x)dx
	\end{aligned}
\end{equation*}
where the last equality follows from $\phi_p(a+\delta-1)=\phi_p(a-\delta-1)=0$ by our construction of $C(x)$.

Thus, if we denote the Lebesgue measure on $\mathbb{T}^2$ by $m$, then for each $p\in[0,1]$, the Markov partition $\mathcal{R}_p$ defined above consists of two rectangles $R_p^0$ with $m(R_p^0)=\frac{5-\sqrt{5}}{10}-\text{Area}_p$ and $R_p^1$ with $m(R_p^1)=\frac{5+\sqrt{5}}{10}+\text{Area}_p$. Since $\text{Area}_p$ is strictly increasing when $p$ increases from $0$ to $1$, $m(R_p^1)$ is strictly increasing when $p$ increases from $0$ to $1$.

In conclusion, for any $p\in[0,1]$, we construct a Markov partition $\mathcal{R}_p=\{R_p^0, R_p^1\}$ for $f_p$, and $m(R_p^1)$ is strictly increasing when $p$ increases from $0$ to $1$.

\

\textbf{Step 3:} Our next step is to prove our main theorem.

\

We define a measurable set $E$ such that
\begin{equation*}
	E=\left\{\begin{array}{l|l}
		(p,z)\in[0,1]\times \mathbb{T}^2 & \begin{array}{l}
			\lim\limits_{n\to\infty}\frac{1}{n}\sum\limits_{k=0}^{n-1}\chi_{R_p^1}(f_p^k(z)) \text{ exists and equals }	m(R_p^1),\\	
			\text{where } \chi_{R_p^1} \text{ is the characteristic function of } R_p^1. 	
		\end{array}
	\end{array}\right\}.
\end{equation*}

Then $E$ is a well-defined and full measure set in $[0,1]\times \mathbb{T}^2$. Because for each fixed $p$, $f_p$ is a smooth area-preserving Anosov diffeomorphism, and hence it is ergodic. By Birkhoff Ergodic theorem, for $m$-a.e.$z\in \mathbb{T}^2$,
\begin{equation*}
	\lim_{n\to\infty}\frac{1}{n}\sum_{k=0}^{n-1}\chi_{R_p^1}(f_p^k(z))=m(R_p^1),
\end{equation*}

Therefore, let $C_p$ denote the torus $\{p\}\times \mathbb{T}^2 \subset [0,1]\times \mathbb{T}^2$. Then the intersection of $E$ with each $C_p$ has two-dimensional Lebesgue measure 1. So, it follows from Fubini's Theorem that $E$ has full three-dimensional Lebesgue measure.

And we note that for any point $(p,z)$, we can  associate a {\em symbol sequence}\/ $(\ldots, b_{-1},b_0,b_1,\ldots)$ with $z$ and $f_p$ such that $b_k = \chi_{R_p^1}(f_p^k(z))$. Then the above limit is exactly
\begin{equation*}
	\lim_{n\to\infty}\frac{b_{0}+\cdots+b_n}{n}.
\end{equation*}

Next, we define a family of curves $\{\Gamma_\beta\}_{\beta\in\mathbb{T}^2}$ as follows. When we take $\epsilon_0$ small enough in step 1, then for any $p\in [0,1]$, we can get a unique homeomorphism $h_p$ on $\mathbb{T}^2$ such that $f_p=h_p\circ f_{0}\circ h_p^{-1}$ since Anosov diffeomorphisms are $C^1$-structurally stable and $f_p$ is a smooth path beginning with $f_0$. Moreover, $h_p$ is continuous with $p$ because the topological conjugacy is continuous with the map $f_p$ (see proof of theorem 2.6.1 in \cite{Katok1995}). So, $h_p(\beta)$ is continuous with both variables. Then for any  $\beta\in \mathbb{T}^2$, we define
\begin{equation*}
	\Gamma_\beta := \{(p,h_p(\beta))\in[0,1]\times \mathbb{T}^2 \; \Big| p\in [0,1]\}.
\end{equation*}
Then $\{\Gamma_\beta\}_{\beta\in \mathbb{T}^2}$ is a family of disjoint curves which fill out  $[0,1] \times \mathbb{T}^2$. And each curve $\Gamma_\beta$ intersects the set $E$ at most a single point. Because if not, we say $\Gamma_\beta$ intersects $E$ more than one point. We may assume two different intersection points to be $(p_1,h_{p_1}(\beta))$ and $(p_2,h_{p_2}(\beta))$. Then $p_1 \neq p_2$. Since these two points are in $E$, we have 
\begin{equation*}
	\lim_{n\to\infty}\frac{1}{n}\sum_{k=0}^{n-1}\chi_{R_{p_i}^1}(f_{p_i}^k(h_{p_i}(\beta)))=m(R_{p_i}^1), \quad i=1,2.
\end{equation*}
And we note that the symbol sequence of $h_{p_1}(\beta)$ under $f_{p_1}$ and the symbol sequence of $h_{p_2}(\beta)$ under $f_{p_2}$ are the same. They both equals the symbol sequence of $\beta$ under $f_0$ denoted by  $(\ldots, \beta_{-1},\beta_0,\beta_1,\ldots)$. This is because $h_p$ is the topological conjugacy and $R_p^1 = h_p (R_0^1)$ by our constructions of Markov partitions. 

To see this, by definition 2.3, we can prove that $h_p(W_0^u((0,0))) = W_p^u((0,0))$ and $h_p(W_0^s((0,0))) = W_p^s((0,0))$. 

Explicitly, if $x\in h_p(W_0^u((0,0)))$, then $h_p^{-1}(x) \in W_0^u((0,0))$ which implies
\begin{equation*}
	\mathrm{dist}(f_0^{-n}(h_p^{-1}(x)),(0,0))\to0\text{ as }n\to\infty.
\end{equation*}
So, by $f_p=h_p\circ f_{0}\circ h_p^{-1}$, we have 
\begin{equation*}
	\mathrm{dist}(f_p^{-n}(x)),h_p(0,0))=\mathrm{dist}(f_p^{-n}(x)),(0,0))\to0\text{ as }n\to\infty.
\end{equation*}
Thus, $x \in W_p^u((0,0))$ i.e. we showed $h_p(W_0^s((0,0)))\subset W_p^u((0,0))$. Similarly, we can show that $W_p^u((0,0))\subset h_p(W_0^s((0,0)))$. And the stable manifold can also be proved in this way.

Then, we note that $h_p$ maps the vertices of $R_0^1$ into themselves as vertices of $R_p^1$. If $z$ is a vertex of $R_0^1$, then $z = L_u \cap L_s$ for some sufficiently small segments $L_u \in W_0^u((0,0))$ and $L_s \in W_0^s((0,0))$ (see Figure 1). Then $h_p(z) =  h_p(L_u) \cap h_p(L_s)$. Note that $h_p(L_u) \subset h_p(W_0^u((0,0))) = W_p^u((0,0))$ and $h_p(L_s) \subset h_p(W_0^s((0,0))) = W_p^s((0,0))$, then by our construction of $R_p^1$ (see Figure 2) we can get  $h_p(z)=h_p(L_u) \cap h_p(L_s)=z$ which is also the vertex of $R_p^1$. 

Since $h_p(W_0^u((0,0))) = W_p^u((0,0))$, $h_p(W_0^s((0,0))) = W_p^s((0,0))$ and $h_p$ is homotopic to identity, we have $h_p$ maps the boundary of $R_0^1$ into the boundary of $R_p^1$ with four vertices fixed. Hence, $h_p (R_0^1)=R_p^1$.

\

Thus the above limits for $i=1,2$ are both equals
\begin{equation*}
	\lim_{n\to\infty}\frac{\beta_{0}+\cdots+\beta_n}{n}.
\end{equation*}

On the other hand, we have $m(R_{p_1}^1)\neq m(R_{p_2}^1)$ because $m(R_p^1)$ is strictly increasing when $p$ increases from 0 to 1. And this gives a contradiction.

\end{document}